\documentclass[12pt]{article}
\usepackage[centertags]{amsmath}
\usepackage{graphicx}
\pdfoutput=1

\begin{document}

\title{On Probability Leakage}
\author{William M. Briggs \\300 E. 71st Apt. 3R, New York, NY 10021\\matt@wmbriggs.com}

\maketitle


\newpage

\vskip .1in\noindent \textsc{Abstract:}

The probability leakage of model M with respect to evidence E is defined.  Probability leakage is a kind of model error.  It occurs when M implies that events $y$, which are impossible given E, have positive probability.   Leakage does not imply model falsification.  Models with probability leakage cannot be calibrated empirically.  Regression models, which are ubiquitous in statistical practice, often evince probability leakage.


\section{Introduction}

We take an objectivist Bayesian view of probability, in the school of Jaynes et al. \cite{Jay2003,Fra2001,Ver1990,Wil2010}.  In this view, probabilities are conditional on evidence E.  For example, this implies for some event R, we could have $\Pr(\operatorname{R} | \operatorname{E}_1) \ne \Pr(R | \operatorname{E}_2)$ for two different sets of evidence E.  For example, if R = ``A five shows", then if $\operatorname{E}_1$ = ``This is a six-sided object with just one side labled 'five' which will be tossed" and $\operatorname{E}_2$ is the same as $\operatorname{E}_1$ except substituting ``just three sides", then $\Pr(\operatorname{R} | \operatorname{E}_1) = 1/6$ and $\Pr(\operatorname{R} | \operatorname{E}_2) = 3/6$. In other words, $\Pr(\operatorname{R})$ is undefined without conditioning evidence. 

A model M is sufficient evidence to define a probability, and of course different models can and do give different probabilities to the same events.  If evidence E suggests that some values of an observable are impossible, yet M gives positive probability to these events, M is said to evince probability leakage with respect to E. 

This has implications for the falsifiability philosophy of \cite{Pop1959,Pop1963} which has led a strange existence in statistics, with many misapprehensions appearing;  e.g. \cite{Gil1971,Spi1974,Daw2004}.  This criterion is clarified with respect to probability leakage.  It is difficult to falsify probability models.

We also take a Bayesian predictive stance, a view which says that all parameters are a nuisance: see {\it inter alia} \cite{Gei1993,Joh1996,JohGei1982,LeeJoh1996}.   This approach allows us to investigate how probability leakage bears on calibration.  Calibration, defined below, is how closely a model's predictions ``match" actual event; e.g. \cite{GneRaf2007}.

Since regression is the most-used statistical model, an example is given which shows how badly these models fare even when all standard diagnostic measures indicate good performance.

\section{Probability Leakage}

\subsection{Definition}

Evidence is gathered or assumed which implies that a model M represents uncertainty in some observable $y$, perhaps conditional on explanatory observable variables $x$ and indexed on parameters $\theta$.  Data are collected: let $z = (y_{old},x_{old})$ be its label.  Ordinary (and inordinate) interest settles on (functions of)

\begin{equation}
  \label{posterior}
   p(\theta|z,\operatorname{M}) \propto p(y_{old}|x_{old},\theta,\operatorname{M})p(\theta|x_{old},\operatorname{M}).
\end{equation}

M is understood to include the evidence that gives the prior, $p(\theta|x_{old}, \operatorname{M})$. The difficulty is that statements like (\ref{posterior}) (or functions of it) cannot be verified.  That is, we never know whether (\ref{posterior}) says something useful or nonsensical about the world. Whatever certainty we have in $\theta$ after seeing data (and assuming M true) tells us nothing directly about $y$.  Indeed, even as $p(\theta|z_n,\operatorname{M}) \to \delta(t)$ (a delta function) as our sample increases, we still do not know with certainty the value of future observables $y$ (given $x$, $z$, and M; and supposing that M is itself not a degenerate distribution).   

Because we assume the truth of M, the following result holds:

\begin{equation}
    \label{predictive}
    p(y|x,z,\operatorname{M}) = \sum_{\theta} p(y|x,z,\theta_i,\operatorname{M}) p(\theta_i|z,\operatorname{M}),
\end{equation}

\noindent with the integral replacing the sum if necessary. This result does not just hold, but it is the logical implication of our previous assumptions; that is, if our assumptions are true, then (\ref{predictive}) must be so.  Even if sole interest is in the posterior (\ref{posterior}), equation (\ref{predictive}) is still implied.  That is to say, given our assumptions, (\ref{predictive}) is deduced.  This logical truth has consequences.

The first is the obvious conclusion that (\ref{posterior}) does not describe $z$ (nor does (\ref{predictive})). All that we want to know about $z$ is in $z$ itself. What remains uncertain are unknown (usually future) values of $y$.  The question is how well does M describe uncertainty in these $y$? Equation (\ref{predictive}) is observable; or, rather, it can be turned into statements which are about observables.   It says that, given M, the old data, and perhaps a new value of $x$, the value of $y$ will take certain values with the calculated probabilities.   Suddenly the world of model verification becomes a possibility (and is explored below).

Now suppose we know, via some evidence E, that $y$ cannot take values outside some set or interval, such as $(y_a,y_b)$. This evidence implies $\Pr(y < y_a | \operatorname{E}) =  \Pr(y > y_b | \operatorname{E}) = 0$.   But if for some value of $x$ (or none is $x$ is null), that $\Pr(y < y_a | x, z, \operatorname{M}) > 0$ or that  $\Pr(y > y_b | x, z, \operatorname{M}) > 0$, then we have a probability leakage; at least, with respect to E.   The probabilities from (\ref{predictive}) are still true, but with respect to M, $z$, and $x$.  They are not true with respect to E if there is leakage. 

This probability leakage is error, but only if we accept E as true. Leakage is a number between 0 (the ideal) and 1 (the model has no overlap with the reality described by E). An example of $y$ with known limits is the GPA of a college student, known to be at some institution strictly between 0 and 4.   These limits form our E.

It's best not to express too rigorous a concern about ``leakage sets," however, as the following example shows.  

Suppose $y$ is the air temperature as measured by a digital thermometer capable of tenth-of-a degree precision. This device, like all tangible devices, will have an upper and lower limit.  It will also because of its resolution give measurements belonging to a finite, discrete set.  This information forms our E.  So not only will the probability (conditional on E) for extreme events be 0, it will also be 0 for all those measurements which cannot register on the device (such as the gaps between the tenths of a degree).  

Now, if we were to use any continuous probability model, such as the normal, to represent uncertainty in this $y$, we would find that the probability for the measurable values (actually observable on the device) to be 0.  The probability leakage of the model with respect to E is complete; it is 1, as high as is possible. Since all real-world measurements known to us are discrete and finite, yet we so often use continuous probability models to represent our uncertainty in these observables, we must turn a blind eye towards leakage of the kind just noted. 

To be clear, whenever E implies $y$ is discrete, yet M gives a continuous representation of $y$, the probability leakage is 1.   Of course,  E could be modified so that it admits of the continuity of $y$, but is it still true that the probability of actual observable events (in real life applications, with respect to any continuous M) is 0. This, of course, is an age-old problem.

So to be interesting we must suppose that $y$, if quantified by a continuous probability model, also lives on the continuum.  We can still acknowledge (in our E) upper and lower limits to $y$, as in the GPA example, where probability leakage is easy to calculate.

\subsection{Model Falsification}

The term falsified is often tossed about, but in a strange and loose way.  A rigorous definition is this: if a M says that a certain event cannot happen, then if that event happens the model M is falsified.  That is, to be falsified M must say that an event is impossible: not unlikely, but impossible. If M implies that some event is merely unlikely, no matter how small this probability, if this event obtains M is not falsified. If M implies that the probability of some event is $\epsilon>0$ then if this event happens, M is not falsified period.  There is thus no escape in the phrase ``Practically impossible," which has the same epistemic properties as ``practically a virgin." See \cite{Wag2004} for how the former phrase can be turned into mathematics.  

Probability leakage does not necessarily falsify M.   If there is incomplete probability leakage, M says certain events have probabilities greater than 0, events which E has says are impossible (have probabilities equal to 0).  If E is true, as we assume it is, then the events M said are possible cannot happen.  But to have falsification of M, we need the opposite: M had to say that events which obtained were impossible.  

Falsification is when events which M says have probability 0 obtain.  This will be the case each time M is a continuous probability distribution and $y$ are real-world, i.e. measured, events; for we have seen that the probability of $y$ taking any measurable value when M is continuous is 0, yet (``real world") E tells us that the probabilities of the $y$ in some discrete set are all greater than 0 (eventually $y$ takes some measurable value).  

Replace the continuous M with a discrete M and model falsification becomes difficult.  Suppose M is a Poisson distribution and $y$ is a count with a known (via some E) upper limit.  M says that the probability of all counts greater than this upper limit have probabilities greater than 0.  But since (if E is true) none of these $y$ will ever obtain, then M will never be falsified.

Box gave us an aphorism which has been corrupted to (in the oral tradition; see \cite{DenKin2006} for a print version), ``All models are wrong."   We can see that this is false: all models are not wrong, i.e. not all are falsified.   They are only wrong, i.e. falsified, if they have complete probability leakage.

\subsection{Calibration}

Calibration  \cite{GneRaf2007} has three components.  Let the empirical distribution of $y$ be represented as $Q$, a distribution we assume is deduced from E, and let the predictive distribution (\ref{predictive}) implied by M be $P$.  In short-hand notation, calibration is defined with respect to  probability
$$
\frac{1}{n} \sum_{i=1}^n Q_i \circ P^{-1}_i(p) \to p;
$$
\noindent exceedance 
$$
\frac{1}{n} \sum Q^{-1}_i \circ P_i(y) \to y;
$$ 
\noindent and marginally
$$
\frac{1}{n} \sum P_i(y) \to \frac{1}{n} \sum Q_i(y), \: \forall y.
$$

If M with respect to E evinces probability leakage, model M cannot be calibrated empirically.   This is easily proved.  In order for  M to be calibrated in probability, the frequency of observed $y$s for which M says that the probability of $y$ or less is $p$ must go to $p$ as the sample increases.   If there is leakage, there will be values of $y$ (call them $y'$) such that $0<P(y')<1$, but which are impossible under E, and which will give $Q(y')=0$ or $Q(y')=1$, making probability calibration impossible.

In order for M to be calibrated in exceedance, each value of $y$ must be as probable under $P$ as $Q$.  If there is leakage, there will be values of $y$ (again $y'$) under E for which $0<P(y')<1$ but $Q^{-1} \circ P(y')= 0$ or 1, making exceedance calibration impossible.

In order for M to be calibrated marginally, its marginal distribution must overlap the empirical distribution.  But if there is probability leakage, there will be impossible values of $y$ (again $y'$) for which $\frac{1}{n} \sum P_i(y) \to p$ for which $\frac{1}{n} \sum Q_i(y') = 0$ or 1, making marginal calibration impossible.

As \cite{GneRaf2007,Sch1989} and others show, if M is to be evaluated by a strictly proper scoring rule, the lack of calibration guarantees that better models than M exist. 

\section{Example}

Statistics as she is practiced---not as we picture her in theoretical perfection---is rife with models exhibiting substantial probability leakage.   This will become obvious from taking just one example, not before published.  It well to point out that the results of most statistical analyses are not published; they remain in private hands and are used in decision making everywhere.

Regression is ubiquitous. The regression model M assumes that $y$ is continuous and that uncertainty in $y$, given some explanatory variables $x$, is quantified by a normal distribution the parameters of which are represented, usually tacitly, by ``flat" (improper) priors.   This M logically implies a T-distribution for (\ref{predictive}), see \cite{BerSmi2000}. This M has the advantage of mimicking standard frequentist results.

I want to emphasize that I do not justify this model for this data; better ones certainly exist. I only claim that regression is often used on data just like this.  The purpose here is only to show how horribly wrong a model can go.

The data is from one of two call center help lines: $y$ is a measure of abandonment, with larger numbers indicating more abandoned calls.  It is impossible that $y$ be less than 0; but it can equal 0.  It has obvious no upper limit. It ranged from 0.08 to 14.1.   It was to be explained by (the $x$s) the number of calls answered, location, and number of absentees (the number of people who were scheduled to work but did not show).  It was expected that greater calls answered would lead to a higher abandonment; it was assumed that the locations differed only in the behavior which would give different abandonments; and it was expected that higher absentees would lead to higher abandonment.   Calls answered could not be less than 0: they ranged from 110 to 2,995.  Absentees could not be less than 0: they ranged from 1 to 14.  There were 52 samples from each location (total 104).

Frequentist model diagnostics gave p-values of 1e-5 for absentees, 0.004 for calls answered, and 1e-5 for the location difference.  The point estimates were of the size and in the directions expected.  Bayesian posteriors on the model parameters showed that each was safely different than 0 (with probabilities at least 0.999).  Visual examination of the ``residuals" showed nothing untoward.   All in all, a very standard regression which performed to expectations and which resulted in a satisfied client. 

But the model is poor for all that because of massive probability leakage. Our E tells us that $y$ cannot be less than 0.  But employing equation (\ref{predictive}) with new values of calls answered and absentees set at 1,200 and 5 (the sample medians) respectively gives Figure 1.    The logical implication of M is that, for these values of $x$, there is about a 38\% chance for values of $y$ less than 0 at Location A.  Even for Location B, there is still a small chance (about 2\%) for values of $y$ less than 0.   

If instead we took the minimum observed values of each $x$  (not pictured), then the probability leakage for Location A is a whopping 92\%, and for Location B it is 50\%.   These are unacceptably large errors. 

In any problem where $x$ is null we can, after observing $z$, calculate the amount of probability leakage. Otherwise, the amount, for any given M, E and $z$, depends on what values of $x$ are expected. 

For example, if we change the M above to a null $x$---i.e. the model assumes only that our uncertainty in $y$ is characterized by a normal distribution---then the probability leakage is the probability, given M and $z$, of $y$ less than 0; which in this case is just north of 10\%. A substantial error. 

But if we keep the $x$s as before, the probability leakage changes with the values of $x$.  This opens the possibility of modeling our uncertainty in probability leakage for a given M, $y$, $x$, and $z$.  It requires, of course, a new model for the leakage as function of the $x$, with all that that entails.  This is left for further research. For now, we only rely on the supposition that the range of $x$s we have seen before are possible, even likely, to realize in new data.

\section{Conclusion}

Probability leakage may be difficult or impossible to compute when limits aren't pertinent, yet there still may be E such that the probability of certain values of $y$ judged to be extremely unlikely are not interesting. An example is when $y$ represents a financial gain or loss.  No obvious limits exist.  We assume $y$'s continuity (via M and E).  Yet we might, given some evidence of the market, say $y$ is so unlikely outside of certain limits that these limits are ``practically impossible."  We have already seen that this phrase is troublesome, but if we make employ it our situation amounts to assuming strict limits as before. 

We might be able to soften the concept of leakage in the absence of limits. An expert might, via probability elicitation, sketch a density (or distribution) for $y$ (possibly given some $x$).  If the distribution (\ref{predictive}) implied by M is ``far" from this elicitation, then either the expert is mistaken or M is.   What is far must be quantified: a possible candidate (among many) is the Kullback-Leibler distance; see for example  \cite{Joh1996}.  This possibility is not explored here.

An objection to the predictive approach is that interest is solely in the posterior (\ref{posterior}); in whether, say, the hypothesis ($\operatorname{H}_{\theta}$) that absentees had an effect on abandonment.  But notice that the posterior does not say with what probability absentees had an effect: it instead says {\it if} M is true and given $z$, the probability that the parameter associated with absentees is this-and-such.  If M is not true (it is falsified), then the posterior has no bearing on $\operatorname{H}_{\theta}$.  In any case, the posterior does not give us  $\Pr(\operatorname{H}_{\theta}|z)$,  it gives $\Pr(\operatorname{H}_{\theta}|\operatorname{M},z)$.  We cannot answer whether $\operatorname{H}_{\theta}$ is likely without referencing M, and M implies (\ref{predictive}). 

Probability leakage is far from the last word in model validation.  It does not answer many questions about the usefulness of models. Nor can it always in isolation tell whether a model is likely true or false: see the book by \cite{Wil2006} for an overview of model validation.  Leakage can, however, give a strong indication of model worth.

\bibliographystyle{plain}
\bibliography{/home/briggs/projects/writing/logic.bib}

\newpage

\begin{figure}[ht]
    \begin{center}
        \scalebox{.8}{\includegraphics{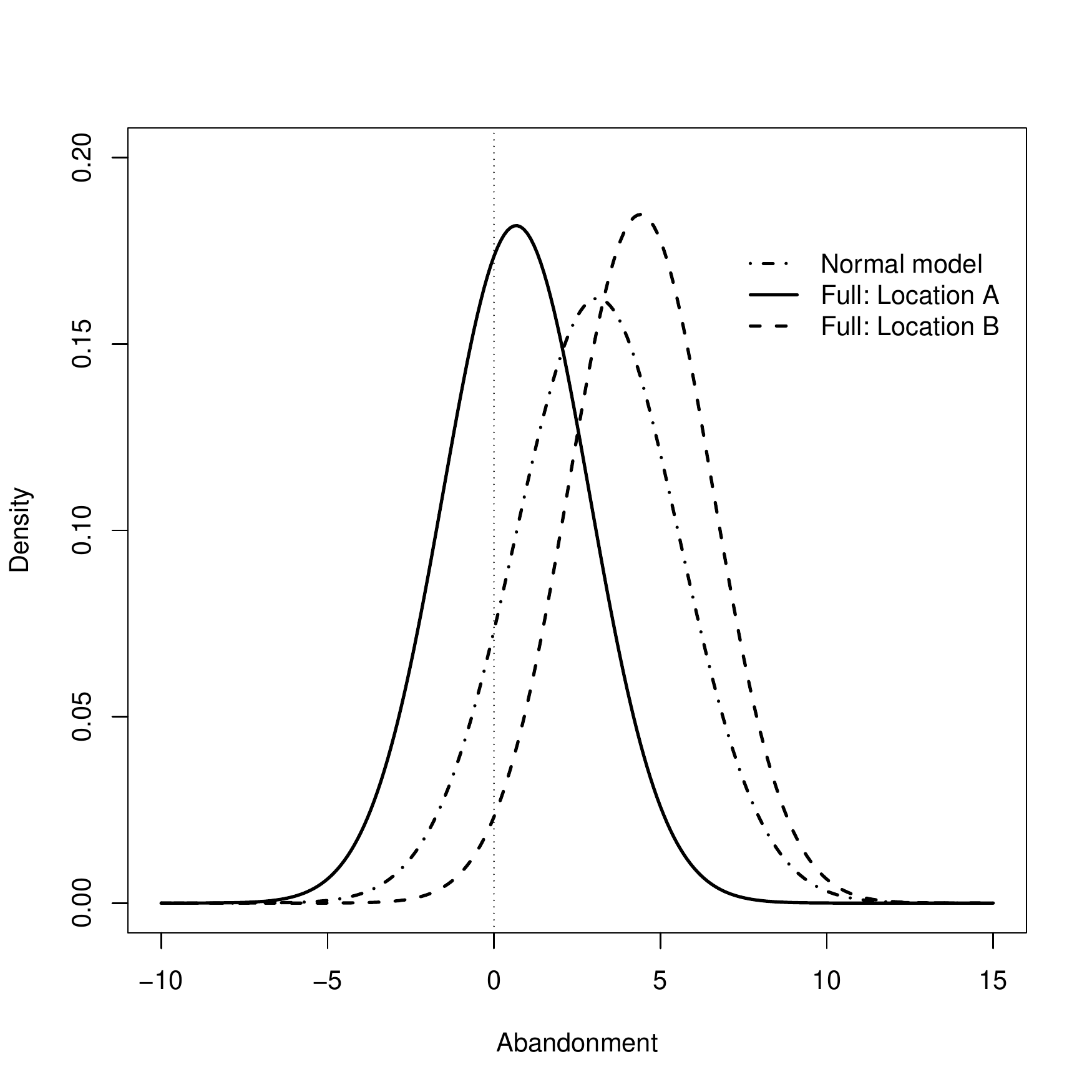}}
    \end{center}
    \caption{The posterior prediction distribution for a normal model with no regressors (dashed-dotted line); for the full regression model at Location A (solid line) and for Location B (dashed line) with the regressors set at their observed medians. The vertical line indicates abandonments of 0, values below which are impossible with respect to E.  Substantial probability is given to impossible values under M. }
\end{figure}

\end{document}